# The Dirichlet Problem for the Minimal Surface System in Arbitrary Codimension

Mu-Tao Wang

## Abstract

Let $\Omega$ be a bounded $C^2$ domain in $\mathbb{R}^n$ and $\phi : \partial \Omega \to \mathbb{R}^m$ be a continuous map. The Dirichlet problem for the minimal surface system asks whether there exists a Lipschitz map $f : \Omega \to \mathbb{R}^m$ with $f|_{\partial \Omega} = \phi$ and with the graph of $f$ a minimal submanifold in $\mathbb{R}^{n+m}$. For $m = 1$, the Dirichlet problem was solved more than thirty years ago by Jenkins-Serrin [13] for any mean convex domains and the solutions are all smooth. This paper considers the Dirichlet problem for convex domains in arbitrary codimension $m$. We prove if $\psi : \overline{\Omega} \to \mathbb{R}^m$ satisfies $8n\delta \sup_\Omega |D^2 \psi| + \sqrt{2} \sup_{\partial \Omega} |D\psi| < 1$, then the Dirichlet problem for $\psi|_{\partial \Omega}$ is solvable in smooth maps. Here $\delta$ is the diameter of $\Omega$. Such a condition is necessary in view of an example of Lawson-Osserman [15]. In order to prove this result, we study the associated parabolic system and solve the Cauchy-Dirichlet problem with $\psi$ as initial data.

## 1 Introduction

There has been emerging attention to minimal submanifolds of higher codimension in recent years. Calibrated submanifolds, e.g. special Lagrangians, arise naturally in string theory and mirror symmetry. This paper considers the Dirichlet problem for the minimal surface equation in arbitrary codimension.

The study of non-parametric minimal hypersurfaces has a long and rich history. Since the seventies, the Dirichlet problem is well-understood due to the work of Jenkins-Serrin [13], De Giorgi [4], and Moser [20]. In contrast, very little is known in higher codimension. In their paper "Non-existence, non-uniqueness, and irregularity of solutions to the minimal surface system" [15], Lawson and Osserman show there exists smooth boundary data $\eta = (\eta^1, \cdots, \eta^m) : S^{n-1} \to \mathbb{R}^m$ such that the Dirichlet problem is not solvable in Lipschitz maps for $R\eta = (R\eta^1, \cdots, R\eta^m)$ when $R$ is large.

Indeed, consider $S^3 \subset \mathbb{C}^2$ with standard complex coordinates $(z_1, z_2)$ and let $\eta : S^3 \to S^2 \subset \mathbb{R}^3$ be the Hopf map given by $\eta(z_1, z_2) = (|z_1|^2 - |z_2|^2, 2z_1\bar{z}_2)$. Scale the image of $\eta$ by $R > 0$ in $\mathbb{R}^3$ to obtain a one-parameter family of boundary maps $R\eta : S^3 \to S^2(R) \subset \mathbb{R}^3$. It was shown in [15] that when $R \geq 4.2$, the Dirichlet problem is not solvable in Lipschitz maps. When $R = \frac{\sqrt{5}}{2}$, a Lipschitz solution exists as a non-parametric minimal cone. The latter example also suggests the regularity theory of non-parametric minimal surfaces differs in codimensions: any Lipschitz minimal hypersurface is smooth.

It is thus natural to conjecture that the solvability of the Dirichlet problem is related to the variation of the boundary map. This paper provides classical solutions to the Dirichlet problem for minimal surface systems in arbitrary codimension assuming such conditions.

**Theorem A.** *Let $\Omega$ be a bounded $C^2$ convex domain in $\mathbb{R}^n$ with diameter $\delta$. If $\psi : \overline{\Omega} \to \mathbb{R}^m$ satisfies $8n\delta \sup_\Omega |D^2\psi| + \sqrt{2} \sup_{\partial\Omega} |D\psi| < 1$, the Dirichlet problem for the minimal surface system is solvable for $\psi|_{\partial\Omega}$ in smooth maps.*

For $x \in \overline{\Omega}$, $|D\psi|(x)$ is the norm of the differential $D\psi(x) : \mathbb{R}^n \to \mathbb{R}^m$ defined by

$$|D\psi|(x) = \sup_{|v|=1} |D\psi(x)(v)|.$$

Similarly,

$$|D^2\psi|(x) = \sup_{|v|=1} |D^2\psi(x)(v,v)|.$$

We solve the Cauchy-Dirichlet problem of the corresponding parabolic system with the initial map $\psi$. The associated parabolic system to the minimal surface equation is the mean curvature flow. It consists of a family of immersions $F : \Omega \times [0, T) \to \mathbb{R}^{n+m}$ satisfying $\frac{dF}{dt} = H$, where $H$ is the mean curvature vector of $F(\Omega, t)$. A maximum principle for higher codimension mean curvature flows was discovered in [21] and [22]. The estimate turns out to be crucial to the existence and regularity problem in higher codimensions.

The solutions to the minimal surface system are critical points of the functional

$$\int \sqrt{\det(I + (Df)^T Df)}$$

where $Df$ is considered as a linear transformation $Df : \mathbb{R}^n \to \mathbb{R}^m$ and $(Df)^T$ is the adjoint of $Df$. The mean curvature flow is the gradient flow of this functional.

A major difference between codimension one and higher codimension is the following: the function $\sqrt{1 + |x|^2}$ is convex for $x \in \mathbb{R}^n$, while

$$\sqrt{\det(I + A^T A)}$$

is not convex in the space of $m \times n$ matrices. The convexity condition is important in the variational approach, see for example Theorem 12.10 in [9] for the existence of Lipschitz solution to the Dirichlet problem in codimension one.

We remark that if we require the minimal submanifold to be Lagrangian, the minimal surface system is reduced to a fully nonlinear scalar equation

$$Im(\det(I + \sqrt{-1}\, D^2 u)) = 0$$

where $I$ = identity matrix and $D^2 u = (\frac{\partial^2 u}{\partial x^i \partial x^j})$ is the Hessian matrix of $u$. The Dirichlet problem for prescribed boundary value of $u$ was solved by Caffarelli, Nirenberg and Spruck in [3].

The paper is organized in the following way. We establish the relation between the parametric and non-parametric form of the mean curvature flow in §2. In §3, a boundary gradient estimate is derived using the initial map as a barrier surface. Theorem A is proved in §4.

I am indebted to Professor R. Hamilton, Professor D. H. Phong, Professor R. Schoen and Professor S.-T. Yau for their constant advice, encouragement and support.

## 2 Preliminaries

Let $\Omega \subset \mathbb{R}^n$ be a domain and $\psi : \Omega \to \mathbb{R}^m$ be a vector-valued function. The graph of $\psi$ is then given by the embedding $I \times \psi : \Omega \to \mathbb{R}^n \times \mathbb{R}^m = \mathbb{R}^{n+m}$ where $I$ is the identity map on $\Omega$.

Let $F : \Omega \times [0, T) \to \mathbb{R}^{n+m}$ be a parametric solution to the Dirichlet problem of the mean curvature flow, i.e

$$\begin{cases} \frac{dF}{dt} = H \\ F|_{\partial \Omega} = I \times \psi|_{\partial \Omega}. \end{cases} \quad (2.1)$$

In terms of local coordinates $(x^1, \cdots, x^n)$ on $\Sigma$, the mean curvature flow is the solution

$$F = F^A(x^1, \cdots, x^n, t), \ A = 1, \cdots, n+m$$

to the following system of parabolic equations

$$\frac{dF}{dt} = \Big( \sum_{i,j=1}^n g^{ij} \frac{\partial^2 F}{\partial x^i \partial x^j} \Big)^\perp \quad (2.2)$$

where $g^{ij} = (g_{ij})^{-1}$ is the inverse of the induced metric $g_{ij} = \frac{\partial F}{\partial x^i} \cdot \frac{\partial F}{\partial x^j}$. The notation $(\cdot)^T$ denotes the tangent part of a vector in $\mathbb{R}^{n+m}$, i.e. the projection onto the subspace spanned by $\{\frac{\partial F}{\partial x^i}\}_{i=1\cdots n}$, while $(\cdot)^\perp$ denotes the normal part or the projection onto the orthogonal complement of the tangent space.

Recall the following formula for the Laplace operator of the metric induced by $F(\cdot, t)$. Note the summation convention, repeated indices are summed over, is adopted in the rest of the article.

**Lemma 2.1**

$$\Delta F \equiv \frac{1}{\sqrt{g}} \frac{\partial}{\partial x^i}\Big(\sqrt{g} g^{ij} \frac{\partial F}{\partial x^j}\Big) = \Big(g^{ij} \frac{\partial^2 F}{\partial x^i \partial x^j}\Big)^\perp$$

*where* $g = \det g_{ij}$.

*Proof.*

We calculate

$$\Delta F = \frac{1}{\sqrt{g}} \frac{\partial}{\partial x^i}\Big(\sqrt{g} g^{ij} \frac{\partial F}{\partial x^j}\Big) = \frac{1}{\sqrt{g}} \frac{\partial}{\partial x^i}(\sqrt{g} g^{ij}) \frac{\partial F}{\partial x^j} + g^{ij} \frac{\partial^2 F}{\partial x^i \partial x^j}.$$

Thus the normal part is

$$(\Delta F)^\perp = \Big(g^{ij} \frac{\partial^2 F}{\partial x^i \partial x^j}\Big)^\perp.$$

To prove the lemma, it suffices to show $\Delta F$ is always in the normal direction. We calculate

$$\Delta F \cdot \frac{\partial F}{\partial x^k} = \frac{1}{\sqrt{g}}\frac{\partial}{\partial x^i}(\sqrt{g}g^{ij}\frac{\partial F}{\partial x^j} \cdot \frac{\partial F}{\partial x^k}) - g^{ij}\frac{\partial F}{\partial x^j} \cdot \frac{\partial^2 F}{\partial x^i \partial x^k}.$$

Now

$$g^{ij}\frac{\partial F}{\partial x^j} \cdot \frac{\partial^2 F}{\partial x^i \partial x^k} = \frac{1}{2}g^{ij}\frac{\partial}{\partial x^k}(g_{ij}).$$

By definition, $\frac{\partial F}{\partial x^j} \cdot \frac{\partial F}{\partial x^k} = g_{jk}$, therefore

$$\Delta F \cdot \frac{\partial F}{\partial x^k} = \frac{1}{\sqrt{g}}\frac{\partial}{\partial x^k}(\sqrt{g}) - \frac{1}{2}g^{ij}\frac{\partial}{\partial x^k}(g_{ij}).$$

Apply the formula

$$\frac{\partial g}{\partial x^k} = gg^{ij}\frac{\partial g_{ij}}{\partial x^k},$$

we obtain

$$\Delta F \cdot \frac{\partial F}{\partial x^k} = 0.$$

□

We now derive a relation between parametric and non-parametric solutions to the mean curvature flow equation. In codimension one case, this was derived in [6].

**Proposition 2.1** *Suppose $F$ is a solution to the Dirichlet problem for mean curvature flow (2.1) and that each $F(\Omega, t)$ can be written as a graph over $\Omega \subset \mathbb{R}^n$. Then there exists a family of diffeomorphism $r_t$ of $\Omega$ such that $\widetilde{F}_t = F_t \circ r_t$ is of the form*

$$\widetilde{F}(x^1, \cdots, x^n) = (x^1, \cdots, x^n, f^1, \cdots, f^m)$$

*and $f = (f^1, \cdots, f^m) : \Omega \times [0, T) \to \mathbb{R}^m$ satisfies*

$$\begin{cases} \frac{df^\alpha}{dt} = g^{ij}\frac{\partial^2 f^\alpha}{\partial x^i \partial x^j} & \alpha = 1, \cdots, m \\ f|_{\partial\Omega} = \psi|_{\partial\Omega} \end{cases} \quad (2.3)$$

where $g^{ij} = (g_{ij})^{-1}$ and $g_{ij} = \delta_{ij} + \sum_\beta \frac{\partial f^\beta}{\partial x^i} \frac{\partial f^\beta}{\partial x^j}$. Conversely, if $f = (f^1, \cdots, f^m) : \Omega \times [0, T) \to \mathbb{R}^m$ satisfies (2.3), then $\widetilde{F} = I \times f$ is a solution to

$$(\frac{d}{dt}\widetilde{F}(x,t))^\perp = \widetilde{H}(x,t).$$

*Proof.* Denote the projections $\pi_1 : \mathbb{R}^{n+m} \to \mathbb{R}^n$ and $\pi_2 : \mathbb{R}^{n+m} \to \mathbb{R}^m$. We first solve for $r_t$ given $f_t = \pi_2 \circ \widetilde{F}_t$. Since the image of $F_t$ is contained in the graph of $f_t$, we have

$$f_t(\pi_1 \circ F) = \pi_2 \circ F_t$$

or

$$\pi_2 \circ \widetilde{F}_t \circ \pi_1 \circ F_t = \pi_2 \circ F_t.$$

It follows

$$r_t = (\pi_1 \circ F_t)^{-1}$$

is the required diffeomorphism. Notice that $r_t$ is the identity map on $\partial \Omega$. Since $\widetilde{F}(x, t) = F(r(x, t), t)$ is simply a reparametrization of $F$,

$$\frac{\partial \widetilde{F}}{\partial t}(x, t) = \frac{\partial F}{\partial t}(r(x, t), t) + dF(\frac{dr}{dt})$$

which implies the normal parts of $\frac{\partial \widetilde{F}}{\partial t}(x, t)$ and $\frac{\partial F}{\partial t}(r(x, t), t)$ coincide.

Thus $\widetilde{F}$ satisfies

$$(\frac{d}{dt}\widetilde{F}(x,t))^\perp = \widetilde{H}(x,t)$$

where $\widetilde{H}(x, t)$ is the mean curvature vector of $\widetilde{F}_t(\Omega)$ at $\widetilde{F}(x, t)$ and $\widetilde{H} = (g^{ij}\frac{\partial^2 \widetilde{F}}{\partial x^i \partial x^j})^\perp$.

Since both $\frac{d}{dt}\widetilde{F}$ and $g^{ij}\frac{\partial^2 \widetilde{F}}{\partial x^i \partial x^j}$ are in $\mathbb{R}^m$ and the map $(\cdot)^\perp$ on $\mathbb{R}^m$ is one-to-one, we obtain

$$\frac{d\widetilde{F}}{dt} = g^{ij}\frac{\partial^2 \widetilde{F}}{\partial x^i \partial x^j}.$$

The last $m$ components give the desired equations for $f$. The other direction follows in the same way.

□

We remark that if we solve for $\widetilde{r} : \Omega \times [0,T) \to \Omega$ that satisfies

$$\frac{\partial}{\partial t}\widetilde{r}(x,t) = -(d\widetilde{F}(\widetilde{r}(x,t),t))^{-1}(\frac{\partial \widetilde{F}}{\partial t}(\widetilde{r}(x,t),t))^T$$

where $d\widetilde{F}(\widetilde{r}(x,t),t) : T_{\widetilde{r}(x,t)}\Omega \to T_{\widetilde{F}(\widetilde{r}(x,t),t)}\widetilde{F}(\Omega)$ and let $F(x,t) = \widetilde{F}(\widetilde{r}(x,t),t)$, then

$$\frac{d}{dt}F(x,t) = H(x,t).$$

In the codimension one case, we have

$$g^{ij}\frac{\partial^2 f}{\partial x^i \partial x^j} = \sqrt{g}\frac{\partial}{\partial x^i}(\frac{1}{\sqrt{g}}\frac{\partial f}{\partial x^i}).$$

This follows from the fact that $g^{ij}$ has a very simple expression

$$g^{ij} = (\delta_{ij} + f_i f_j)^{-1} = \delta_{ij} - \frac{f_i f_j}{1 + |Df|^2}.$$

The Cauchy-Dirichlet problem for the equation

$$\frac{df}{dt} = \sqrt{1+|Df|^2}\frac{\partial}{\partial x^i}(\frac{1}{\sqrt{1+|Df|^2}}\frac{\partial f}{\partial x^i})$$

has been studied by Huisken [12] and Lieberman [16].

The following zeroth order estimate is a direct consequence of the maximum principle .

**Proposition 2.2** *Let $f : \Omega \times [0,T) \to \mathbb{R}^m$ be a solution to equation (2.3), if $\sup_{\Omega \times [0,T)}|Df|$ is bounded then*

$$\sup_{\Omega \times [0,T)} f^\alpha \leq \sup_\Omega \psi^\alpha$$

*where $\psi$ is the initial map in equation (2.3).*

# 3 Boundary gradient estimates

The key point in the boundary gradient estimate is to construct barrier functions. The estimate in this section makes use of the initial data

$$\psi : \Omega \to \mathbb{R}^m$$

as a barrier surface. We shall adopt the non-parametric version of the mean curvature flow (2.3) throughout this section.

Let $\Gamma$ be an $(n-1)$ dimensional submanifold of $\mathbb{R}^n$ and $d_\Gamma(\cdot)$ be the distance function to $\Gamma$. Let $f = (f^1, \cdots, f^m)$ be a solution of equation (2.3). For each $\alpha = 1 \cdots m$, consider the following function defined on $\mathbb{R}^n$.

$$S(x^1, \cdots, x^n, t) = \nu \log(1 + k d_\Gamma) - (f^\alpha - \psi^\alpha)$$

where $\nu, k > 0$ are to be determined.

We denote $g^{ij} \frac{\partial^2}{\partial x^i \partial x^j}$ by $\Delta$ in this section where $g^{ij} = (g_{ij})^{-1}$ and $g_{ij} = \delta_{ij} + \sum_\beta \frac{\partial f^\beta}{\partial x^i} \frac{\partial f^\beta}{\partial x^j}$. Direct calculation shows

$$-\Delta \log(1 + k d_\Gamma) = \frac{k}{1 + k d_\Gamma}(-\Delta d_\Gamma) + \frac{k^2}{(1 + k d_\Gamma)^2} g^{ij} \frac{\partial d_\Gamma}{\partial x^i} \frac{\partial d_\Gamma}{\partial x^j}$$

Therefore $S$ satisfies the following evolution equation.

$$(\frac{d}{dt} - \Delta)S = \frac{\nu k}{1 + k d_\Gamma}(-\Delta d_\Gamma) + \frac{\nu k^2}{(1 + k d_\Gamma)^2} g^{ij} \frac{\partial d_\Gamma}{\partial x^i} \frac{\partial d_\Gamma}{\partial x^j} - \Delta \psi^\alpha. \quad (3.1)$$

Now we prove the boundary gradient estimate for convex domains.

**Theorem 3.1** *Let $\Omega$ be a bounded $C^2$ convex domain in $\mathbb{R}^n$. Suppose the flow (2.3) exists smoothly on $\Omega \times [0, T)$, then the following boundary gradient estimate holds:*

$$|Df| < 4n\delta(1 + \xi) \sup_\Omega |D^2 \psi| + \sqrt{2} \sup_{\partial \Omega} |D\psi| \quad on \ \partial \Omega \times [0, T)$$

*where $\xi = \sup_{\Omega \times [0,T)} |Df|^2$ and $\delta$ is the diameter of $\Omega$.*

*Proof.* We take $\Gamma$ to be the supporting $n-1$ dimensional hyperplane $P$ at a boundary point $p$. Since $d_P$ is a linear function, $\Delta d_P = 0$ and equation (3.1) becomes

$$(\frac{d}{dt} - \Delta)S = \frac{\nu k^2}{(1+kd_P)^2}g^{ij}\frac{\partial d_P}{\partial x^i}\frac{\partial d_P}{\partial x^j} - \Delta\psi^\alpha. \tag{3.2}$$

Since the eigenvalues of $g^{ij}$ are between $\frac{1}{1+\xi}$ and 1 and $|Dd_P| = 1$, we obtain

$$g^{ij}\frac{\partial d_P}{\partial x^i}\frac{\partial d_P}{\partial x^j} \geq \frac{1}{1+\xi}.$$

Now

$$\frac{\nu k^2}{(1+kd_P)^2}g^{ij}\frac{\partial d_P}{\partial x^i}\frac{\partial d_P}{\partial x^j} \geq \frac{\nu k^2}{(1+k\delta)^2}\frac{1}{1+\xi}$$

because $d_P(y) \leq |y - p| \leq \delta$ for any $y \in \Omega$. On the other hand,

$$|\Delta\psi^\alpha| = |g^{ij}\frac{\partial^2\psi}{\partial x^i \partial x^j}| \leq n|D^2\psi|.$$

Now we require

$$\frac{\nu k^2}{(1+k\delta)^2}\frac{1}{1+\xi} \geq n\sup_\Omega |D^2\psi|. \tag{3.3}$$

In view of (3.2), the condition (3.3) guarantees $(\frac{d}{dt} - \Delta)S \geq 0$ on $[0, T)$. Notice that on the boundary of $\Omega$, we have $S > 0$ except $S = 0$ at $p$ by convexity. On the other hand, $S \geq 0$ on $\Omega$ at $t = 0$. It follows from the strong maximum principle that $S > 0$ on $\Omega \times (0, T)$. Likewise we can apply this procedure to $S' = \nu\log(1 + kd_\Gamma) + (f^\alpha - \psi^\alpha)$. Therefore at the boundary point $p$, the normal derivative satisfies

$$|\frac{\partial(f^\alpha - \psi^\alpha)}{\partial n}|(p) \leq \lim_{d_P(x) \to 0}\frac{|f^\alpha(x) - \psi^\alpha(x)|}{d_P(x)} < \lim_{d_P(x) \to 0}\frac{\nu\log(1 + kd_P(x))}{d_P(x)} = \nu k.$$

We may assume $\frac{\partial f^\alpha}{\partial n} = 0$ for all $\alpha$ except $\alpha = 1$ by changing coordinates of $\mathbb{R}^m$ to obtain

$$|\frac{\partial f}{\partial n}| < \nu k + |\frac{\partial \psi}{\partial n}|.$$

For $x \in \partial\Omega$, define $|D^{\partial\Omega}f|(x) = \sup_v |Df(x)(v)|$ where the sup is taken over all unit vectors $v$ tangent to $\partial\Omega$. The Dirichlet condition implies

$$|D^{\partial\Omega}f| = |D^{\partial\Omega}\psi|$$

on $\partial\Omega$.

Therefore

$$|Df| < \sqrt{(\nu k + |\frac{\partial \psi}{\partial n}|)^2 + |D^{\partial\Omega}\psi|^2} \leq \nu k + \sqrt{2}|D\psi|$$

on $\partial\Omega$.

Now we can minimize $\nu k$ subject to the constraint (3.3). The minimum is achieved when $k = \delta^{-1}$ and $\nu k = 4n\delta(1+\xi)\sup_\Omega |D^2\psi|$. The theorem is proved.

□

# 4 Proof of Theorem A

*Proof.*

We divide the proof of Theorem A into several steps:

1. We can prove short time existence as Theorem 8.2 in [17]. By the Schauder fixed point theorem, the solvability of equation (2.3) reduces to the estimates of the solution $(f^\alpha)$ to

$$\begin{cases} \frac{d}{dt}f^\alpha = \tilde{g}^{ij}\frac{\partial^2 f^\alpha}{\partial x^i \partial x^j}, \ \alpha = 1, \cdots, m \\ f|_{\partial\Omega} = \psi|_{\partial\Omega} \end{cases} \quad (4.1)$$

where $\tilde{g}^{ij} = (\tilde{g}_{ij})^{-1}$ and $\tilde{g}_{ij} = \delta_{ij} + \sum_\beta \frac{\partial u^\beta}{\partial x^i}\frac{\partial u^\beta}{\partial x^j}$, for any $u = (u^\alpha)$ with uniform $C^{1,\gamma}$ bound.

   Notice that equation (4.1) is a decoupled system of linear parabolic equations. The equation is uniform parabolic and the required estimate follows from linear theory for scalar equations.

2. Denote the graph of $f_t$ by $\Sigma_t$. We show that $|Df_t| < 1$ holds under the assumption of the theorem. Instead of calculating the evolution equation of the differential of $f$, we consider the restriction of following symmetric bilinear form defined on $\mathbb{R}^{n+m}$,

$$P(X,Y) = \langle \pi_1(X), \pi_1(Y) \rangle - \langle \pi_2(X), \pi_2(Y) \rangle$$

where $\pi_1 : \mathbb{R}^{n+m} = \mathbb{R}^n \times \mathbb{R}^m \to \mathbb{R}^n$ and $\pi_2 : \mathbb{R}^{n+m} \to \mathbb{R}^m$ are the projections onto the first and the second factor, respectively. It is clear that $|Df_t| < 1$ if and only if the restriction of $P$ to the tangent space of $\Sigma_t$ is positive definite.

We shall use the parametric version of the mean curvature flow to derive the evolution equation of $P$. Let $F : \Omega \times [0,T) \to \mathbb{R}^{n+m}$ be a mean curvature flow. Fix a coordinate system $\{x^i\}$ on $\Omega$ and denote the restriction of $P$ to $\Sigma_t$ by $P_{ij} = P(\frac{\partial F}{\partial x^i}, \frac{\partial F}{\partial x^j})$. $\tilde{\Delta}$ denotes the rough Laplacian on tensors, so $\tilde{\Delta}P = g^{kl}\nabla_k\nabla_l P$ where $g^{kl} = g_{kl}^{-1}$, $g_{kl} = \frac{\partial F}{\partial x^k} \cdot \frac{\partial F}{\partial x^l}$ is the induced metric and $\nabla_k$ is the covariant derivative on $\Sigma_t$ with respect to $\frac{\partial F}{\partial x^k}$.

We shall look at the evolution equation of $P_{ij}$ at a space-time point $(p,t)$. We choose an orthonormal basis $\{e_\alpha\}_{\alpha=1\cdots m}$ for the normal space at $(p,t)$ and denote the $\alpha$ component of the second fundamental form by $h_{\alpha ki} = \frac{\partial^2 F}{\partial x^k \partial x^i} \cdot e_\alpha$ and of the mean curvature vector by $H_\alpha = H \cdot e_\alpha$. Then $P_{ij}$ satisfies

$$\begin{aligned}(\frac{d}{dt} - \tilde{\Delta})P_{ij} &= g^{nm}(-h_{\alpha in}H_\alpha + g^{kl}h_{\alpha ki}h_{\alpha ln})P_{mj} \\ &+ g^{nm}(-h_{\alpha jn}H_\alpha + g^{kl}h_{\alpha kj}h_{\alpha ln})P_{im} \\ &- 2g^{kl}h_{\alpha ki}h_{\beta lj}P_{\alpha\beta}\end{aligned} \quad (4.2)$$

where $P_{\alpha\beta} = P(e_\alpha, e_\beta)$.

Equation (4.2) is essentially derived in [21] (see equation (2.3)) and we derive it again in the appendix for completeness.

We claim if we assume

$$8n\delta \sup_\Omega |D^2\psi| + \sqrt{2}\sup_{\partial\Omega}|D\psi| < 1$$

then $\sup_\Omega |Df_t| < 1$ as long as the flow exists smoothly. By integration along a path in $\Omega$, we have $\sup_\Omega |Df_0| = \sup_\Omega |D\psi| < 1$ initially. Suppose $|Df_t| = 1$ for the first time at $t_0$ at some $p_0$ in some tangent direction $v$. Theorem 3.1 implies $\sup_{\partial\Omega}|Df_t| < 1$. Therefore $p_0$ must be an interior point and $v$ is a null-vector of $P_{ij}$. We then apply Hamilton's maximum principle (Theorem 9.1 in [10], see also [11]) for tensors to equation (4.2) and show $P_{ij}$ being positive definite is preserved. Denote the right hand side of (4.2) by $N_{ij}$. We need to show whenever $P_{ij}v^i = 0$ for all $j$, we have $N_{ij}v^iv^j \geq 0$.

It is easy to see the terms

$$[g^{nm}(-h_{\alpha in}H_\alpha + g^{kl}h_{\alpha ki}h_{\alpha ln})P_{mj} + g^{nm}(-h_{\alpha jn}H_\alpha + g^{kl}h_{\alpha kj}h_{\alpha ln})P_{im}]v^iv^j = 0$$

vanish for such $v^i$. The key point now is when $P_{ij}$ is non-negative definite, $P_{\alpha\beta}$ is non-positive definite. This can be seen by applying singular value decomposition to the linear transformation $Df : \mathbb{R}^n \to \mathbb{R}^m$. For simplicity we assume $n = m$, the general case follows in the same way except the notation is more complicated. At any point on $\Sigma_t$, let $\{\lambda_i\}_{i=1\cdots n}$ be the singular values of $Df$, then there exist an orthonormal basis $\{a_i\}_{i=1\cdots n}$ for $\mathbb{R}^n$ and an orthonormal basis $\{a_{n+i}\}_{i=1\cdots n}$ for $\mathbb{R}^m$ such that
$$Df(a_i) = \lambda_i a_{n+i}.$$
Then $\{e_i = \frac{1}{\sqrt{1+\lambda_i^2}}(a_i + \lambda_i a_{n+i})\}$ and $\{e_{n+i} = \frac{1}{\sqrt{1+\lambda_i^2}}(a_{n+i} - \lambda_i a_i)\}$ form orthonormal bases of the tangent and normal space, respectively. In these bases, the restriction of $P$ to the tangent space is $\frac{1-\lambda_i^2}{1+\lambda_i^2}\delta_{ij}$ and the restriction to the normal space is $\frac{\lambda_i^2-1}{1+\lambda_i^2}\delta_{ij}$. The same method can be used to proved $P_{ij} - \epsilon g_{ij} > 0$, or equivalently $|Df|^2 < \frac{1-\epsilon}{1+\epsilon}$ for any $\epsilon > 0$ is preserved along the flow.

3. We can prove long-time existence by blow-up analysis as in [21] and [22]. We first recall a formula for mean curvature flows from [22] and [24]. Let $\Omega_1$ be the volume form of $\mathbb{R}^n$ and

$$*\Omega_1 = \frac{\Omega_1(\frac{\partial F}{\partial x^1}, \cdots, \frac{\partial F}{\partial x^n})}{\sqrt{\det g_{ij}}}$$

be the Jacobian of the projection $\pi_1$ restricted to $\Sigma_t$. $*\Omega_1 = \frac{1}{\sqrt{\prod(1+\lambda_i^2)}}$ in terms of the singular values of $Df$. The evolution equation for $*\Omega_1$ in derived in [22] (equation (3.8)). To obtain the equation for $\ln *\Omega_1$ from equation (3.8) of [22], simply follow the calculation in Proposition 2.1 of [24]. We obtain

$$(\frac{d}{dt}-\Delta)(\ln *\Omega_1) = -\{\sum_{\alpha,l,k} h_{\alpha lk}^2 + \sum_{k,i}\lambda_i^2 h_{n+i,ik}^2 + 2\sum_{k,i<j}\lambda_i\lambda_j h_{n+i,jk}h_{n+j,ik}\}. \quad (4.3)$$

Now $|Df| < 1$ implies $|\lambda_i \lambda_j| < 1$ and the right hand side can be completed square to

$$(\frac{d}{dt}-\Delta)(\ln *\Omega_1) \leq -\epsilon_1 |A|^2$$

for some $\epsilon_1 > 0$. The argument in Proposition 6.1 of [21] (see also Theorem A of [22]) shows $|A|^2$ vanishes on a parabolic blow-up limit, so White's regularity theorem [25] gives $C^{2,\gamma}$ bound for some $\gamma$. All the higher order estimates can be derived as in the codimension one case, see for example [6]. In conclusion, we obtain finite time regularity $f_t \in C^\infty(\Omega) \cap C^1(\overline{\Omega})$. Therefore, the solution to the mean curvature flow (2.3) exists smoothly in $[0, \infty)$.

4. The area element $\sqrt{g}$ along the mean curvature flow satisfies $\frac{d}{dt}\sqrt{g} = -|H|^2\sqrt{g}$. Integrating over space and time we find that

$$\int_0^\infty \int_{\Sigma_t} |H|^2$$

is bounded. We can find a sequence $\Sigma_{t_i}$ such that $\int_{\Sigma_{t_i}} |H|^2 \to 0$. Since we have gradient bound too, we may extract a subsequence $t_i$ such that $\Sigma_{t_i}$ converges to a Lipschitz graph with $\int |H|^2 = 0$ and $|Df| < 1$.

5. Interior Regularity of the limit is proved in the next theorem. Boundary regularity follows from Allard's Theorem [2], see Theorem 2.3 in [15].

$\square$

Lawson-Osserman's example of non-parametric minimal cone also demonstrates the difference between codimension one and higher codimension in regularity theory. In codimension one, a Lipschitz non-parametric minimal cone is smooth. A scalar solution to a quasilinear uniform elliptic equation is smooth, while a solution to a uniform elliptic system may not be. We do have the following regularity theorem when $|Df|$ is better controlled.

**Theorem 4.1** *Let $\Omega$ be a bounded domain in $\mathbb{R}^n$ and $f : \Omega \to \mathbb{R}^m$ be a Lipschitz solution to the minimal surface system such that $|\lambda_i \lambda_j| \leq 1 - \epsilon_2$ for any $i \neq j$ for some $\epsilon_2 > 0$, then $f$ is smooth.*

*Proof.* This follows from the Bernstein type theorem proved in [24]. We show the tangent cone at any possible singular point $p_0$ is flat. Indeed, any blow up limit at $p_0$ is a minimal cone by the monotonicity formula. We then apply dimension reduction to show the only possible singularity is the vertex of the cone. Notice that the condition $|\lambda_i \lambda_j| \leq 1 - \epsilon$ is preserved under blow-up and thus remains true on the minimal cone. The proof of Theorem 1.1 in [24] (using essentially the stationary form of (4.3)) applies to this situation, thus the tangent cone at $p_0$ is flat and then Allard's regularity theorem [1] implies $p_0$ is a regular point. $\square$

## 5 Appendix

Let $F : \Omega \times [0, T) \to \mathbb{R}^{n+m}$ be a parametrized mean curvature flow and let $\Sigma_t = F(\Omega, t)$. Given a bilinear form $P = P_{AB}$ defined on $\mathbb{R}^{n+m}$. The restriction of $P$ to $\Sigma_t$ is denoted by $P_{ij} = P_{AB} \frac{\partial F^A}{\partial x^i} \frac{\partial F^B}{\partial x^j}$. $P_{ij}$ is a time-dependent two-tensor defined on $\Sigma_t$ and we consider the rough Laplacian $\tilde{\Delta} P = g^{kl} \nabla_k \nabla_l P$ with respect to the time-dependent induced metric.

First we have

$$\tilde{\Delta} P(\frac{\partial F}{\partial x^i}, \frac{\partial F}{\partial x^j}) = g^{kl} \nabla_k \nabla_l P(\frac{\partial F}{\partial x^i}, \frac{\partial F}{\partial x^j})$$
$$= g^{kl}[\frac{\partial}{\partial x^k}(\nabla_l P(\frac{\partial F}{\partial x^i}, \frac{\partial F}{\partial x^j})) - \nabla_l P(\nabla_k \frac{\partial F}{\partial x^i}, \frac{\partial F}{\partial x^j}) - \nabla_l P(\frac{\partial F}{\partial x^i}, \nabla_k \frac{\partial F}{\partial x^j})]$$

by the Leibnitz rule of the covariant derivative $\nabla$ on $\Sigma_t$. To simply the calculation, we can choose our coordinate so that the connection term $\nabla_k \frac{\partial F}{\partial x^i} = (\frac{\partial^2 F}{\partial x^k \partial x^i})^T = 0$ at a space-time point $(p, t)$.

By Leibnitz rule again,

$$\tilde{\Delta} P(\frac{\partial F}{\partial x^i}, \frac{\partial F}{\partial x^j}) = g^{kl} \frac{\partial}{\partial x^k}[\frac{\partial}{\partial x^l} P(\frac{\partial F}{\partial x^i}, \frac{\partial F}{\partial x^j}) - P(\nabla_l \frac{\partial F}{\partial x^i}, \frac{\partial F}{\partial x^j}) - P(\frac{\partial F}{\partial x^i}, \nabla_l \frac{\partial F}{\partial x^j})] \tag{5.1}$$

We compute the terms in the bracket,

$$\frac{\partial}{\partial x^l} P(\frac{\partial F}{\partial x^i}, \frac{\partial F}{\partial x^j}) - P(\nabla_l \frac{\partial F}{\partial x^i}, \frac{\partial F}{\partial x^j}) - P(\frac{\partial F}{\partial x^i}, \nabla_l \frac{\partial F}{\partial x^j})$$
$$= P(\frac{\partial^2 F}{\partial x^l \partial x^i}, \frac{\partial F}{\partial x^j}) + P(\frac{\partial F}{\partial x^i}, \frac{\partial^2 F}{\partial x^l \partial x^j}) - P(\nabla_l \frac{\partial F}{\partial x^i}, \frac{\partial F}{\partial x^j}) - P(\frac{\partial F}{\partial x^i}, \nabla_l \frac{\partial F}{\partial x^j})$$
$$= P((\frac{\partial^2 F}{\partial x^l \partial x^i})^\perp, \frac{\partial F}{\partial x^j}) + P(\frac{\partial F}{\partial x^i}, (\frac{\partial^2 F}{\partial x^l \partial x^j})^\perp).$$

Plug this into equation (5.1), we continue

$$\frac{\partial}{\partial x^k}[P((\frac{\partial^2 F}{\partial x^l \partial x^i})^\perp, \frac{\partial F}{\partial x^j})]$$
$$= P(\frac{\partial}{\partial x^k}(\frac{\partial^2 F}{\partial x^l \partial x^i})^\perp, \frac{\partial F}{\partial x^j}) + P((\frac{\partial^2 F}{\partial x^l \partial x^i})^\perp, (\frac{\partial^2 F}{\partial x^k \partial x^j})^\perp)$$
$$= P((\frac{\partial}{\partial x^k}(\frac{\partial^2 F}{\partial x^l \partial x^i})^\perp)^T + (\frac{\partial}{\partial x^k}(\frac{\partial^2 F}{\partial x^l \partial x^i})^\perp)^\perp, \frac{\partial F}{\partial x^j}) + P((\frac{\partial^2 F}{\partial x^l \partial x^i})^\perp, (\frac{\partial^2 F}{\partial x^k \partial x^j})^\perp)$$

We can commute the derivatives as the following

$$P((\frac{\partial}{\partial x^k}(\frac{\partial^2 F}{\partial x^l \partial x^i})^\perp)^\perp, \frac{\partial F}{\partial x^j})$$
$$= P((\frac{\partial^3 F}{\partial x^k \partial x^l \partial x^i})^\perp)^\perp - (\frac{\partial}{\partial x^k}\nabla_l \frac{\partial F}{\partial x^i})^\perp, \frac{\partial F}{\partial x^j})$$
$$= P((\frac{\partial}{\partial x^i}(\frac{\partial^2 F}{\partial x^k \partial x^l})^\perp)^\perp + (\frac{\partial}{\partial x^i}\nabla_k \frac{\partial F}{\partial x^l})^\perp - (\frac{\partial}{\partial x^k}\nabla_l \frac{\partial F}{\partial x^i})^\perp, \frac{\partial F}{\partial x^j}).$$

Now
$$(\frac{\partial}{\partial x^i}\nabla_k \frac{\partial F}{\partial x^l})^\perp = (\frac{\partial}{\partial x^k}\nabla_l \frac{\partial F}{\partial x^i})^\perp = 0$$

at $p$ since $\nabla_k \frac{\partial F}{\partial x^l}$ is a tangent vector and $g^{kl}(\frac{\partial}{\partial x^i}(\frac{\partial^2 F}{\partial x^k \partial x^l})^\perp)^\perp = (\frac{\partial H}{\partial x^i})^\perp$.

Therefore

$$\tilde{\Delta} P(\frac{\partial F}{\partial x^i}, \frac{\partial F}{\partial x^j}) = P((\frac{\partial H}{\partial x^i})^\perp, \frac{\partial F}{\partial x^j}) + P(\frac{\partial F}{\partial x^i}, (\frac{\partial H}{\partial x^j})^\perp)$$
$$+ g^{kl}\{P((\frac{\partial}{\partial x^k}(\frac{\partial^2 F}{\partial x^l \partial x^i})^\perp)^T, \frac{\partial F}{\partial x^j}) + P(\frac{\partial F}{\partial x^i}, (\frac{\partial}{\partial x^k}(\frac{\partial^2 F}{\partial x^l \partial x^j})^\perp)^T)$$
$$+ P((\frac{\partial^2 F}{\partial x^l \partial x^i})^\perp, (\frac{\partial^2 F}{\partial x^k \partial x^j})^\perp) + P((\frac{\partial^2 F}{\partial x^k \partial x^i})^\perp, (\frac{\partial^2 F}{\partial x^l \partial x^j})^\perp)\}$$

On the other hand,

$$\frac{d}{dt} P(\frac{\partial F}{\partial x^i}, \frac{\partial F}{\partial x^j}) = P(\frac{\partial H}{\partial x^i}, \frac{\partial F}{\partial x^j}) + P(\frac{\partial F}{\partial x^i}, \frac{\partial H}{\partial x^j})$$

Therefore

$$(\frac{d}{dt} - \tilde{\Delta}) P(\frac{\partial F}{\partial x^i}, \frac{\partial F}{\partial x^j}) = P((\frac{\partial H}{\partial x^i})^T, \frac{\partial F}{\partial x^j}) + P(\frac{\partial F}{\partial x^i}, (\frac{\partial H}{\partial x^j})^T)$$
$$- g^{kl}\{P((\frac{\partial}{\partial x^k}(\frac{\partial^2 F}{\partial x^l \partial x^i})^\perp)^T, \frac{\partial F}{\partial x^j}) + P(\frac{\partial F}{\partial x^i}, (\frac{\partial}{\partial x^k}(\frac{\partial^2 F}{\partial x^l \partial x^j})^\perp)^T)$$
$$+ P((\frac{\partial^2 F}{\partial x^l \partial x^i})^\perp, (\frac{\partial^2 F}{\partial x^k \partial x^j})^\perp) + P((\frac{\partial^2 F}{\partial x^k \partial x^i})^\perp, (\frac{\partial^2 F}{\partial x^l \partial x^j})^\perp)\}$$

To put this in the form of equation (4.2), we simply write out each term in $H_\alpha$ and $h_{\alpha ij}$. For example,

$$(\frac{\partial H}{\partial x^i})^T = \langle \frac{\partial H}{\partial x^i}, \frac{\partial F}{\partial x^n} \rangle g^{nm} \frac{\partial F}{\partial x^m} = -\langle H, \frac{\partial^2 F}{\partial x^n \partial x^i} \rangle g^{nm} \frac{\partial F}{\partial x^m} = -H_\alpha h_{\alpha in} g^{nm} \frac{\partial F}{\partial x^m}.$$